\documentclass[%
 reprint,
 amsmath,amssymb,
 aps,
]{revtex4-1}
\usepackage{mathtools}
\usepackage{graphicx}
\usepackage{dcolumn}
\usepackage{bm}
\usepackage{hyperref}

\begin{document}

\title{Atomic subgraphs and the statistical mechanics of networks}

\author{Anatol E. Wegner}

 \email{a.wegner@ucl.ac.uk}
\affiliation{%
 Department of Statistical Science,
 University College London
}%

\author{Sophia Olhede}

\affiliation{
 Institute of Mathematics, Chair of Statistical Data Science, EPFL
}%
\affiliation{
 Department of Statistical Science,
 University College London
}%
\email{sofia.olhede@epfl.ch}

\date{\today}

\begin{abstract}

 We develop random graph models where graphs are generated by connecting not only pairs of vertices by edges but also larger subsets of vertices by copies of small atomic subgraphs of arbitrary topology. This allows the for the generation of graphs with extensive numbers of triangles and other network motifs commonly observed in many real world networks. More specifically we focus on maximum entropy ensembles under constraints placed on the counts and distributions of atomic subgraphs and derive general expressions for the entropy of such models. We also present a procedure for combining distributions of multiple atomic subgraphs that enables the construction of models with fewer parameters. Expanding the model to include atoms with edge and vertex labels we obtain a general class of models that can be parametrized in terms of basic building blocks and their distributions that includes many widely used models as special cases. These models include random graphs with arbitrary distributions of subgraphs, random hypergraphs, bipartite models, stochastic block models, models of multilayer networks and their degree corrected and directed versions. We show that the entropy for all these models can be derived from a single expression that is characterized by the symmetry groups of atomic subgraphs.

\end{abstract}

\maketitle

\section{Introduction}
Random graph models are fundamental to understanding the interplay between topological features of networks and the effect of topological features on dynamical processes on graphs. Traditionally, random graph models have concentrated on specific features commonly observed in real world networks such as community structure \cite{Holland1983StochasticSteps}, heterogeneous degree distributions \cite{Newman2001RandomApplications} and short geodesic path lengths \cite{Watts1998Collectivesmall-worldnetworks}. Although this resulted in a large variety of models, most are limited in their scope and only aim to replicate a small subset of features while being unrealistic with respect to others. More recently models that combine multiple features have become more prominent in the field. For instance the degree corrected stochastic block model (DC-SBM) \cite{Karrer2011StochasticNetworks,Peixoto2012EntropyEnsembles}, which unifies the stochastic block model (SBM) for community structures and the configuration model for networks with heterogeneous degree distributions, is a much better fit for many empirical networks than the SBM \cite{Peixoto2013ParsimoniousNetworks,Peixoto2014HierarchicalNetworks}. The DC-SBM has been be further generalized various network types \cite{Peixoto2015ModelGroups,Peixoto2015InferringNetworks,Peixoto2018ReconstructingErrors,Peixoto2019NetworkDynamics} and has produced a general framework for the statistical inference of network communities as well as methods for discriminating between alternative representations of networks via models selection in a principled and consistent manner \cite{Peixoto2015ModelGroups,Peixoto2015InferringNetworks}. However, modelling the prevalence of triangles and other motifs \cite{Milo2002NetworkNetworks} observed in many real world networks still is a major challenge in developing realistic random graph models. 

In this article we seek to formulate a class of analytically tractable models that not only can generate graphs with realistic subgraph structures but also provides a unified description of a large variety of models each aimed at modelling seemingly unrelated features of networks by showing that they can be described in terms of atomic building blocks and constraints placed on their distributions. The statistical ensembles we obtain are in many instances special cases or generalizations of previously proposed models which can be solved analytically for many of their properties including topological phase transitions, subgraph distributions and percolation properties \cite{Bollobas2011SparseClustering,Karrer2010RandomSubgraphs}. We focus our efforts on the entropy and likelihood due to their relevance to statistical inference. In doing so we seek to provide a general class of models that can be used to infer statistically significant features of networks and for discriminating between alternative representations of networks via model selection.

Maximum entropy models \cite{Cimini2019TheNetworks} offer a general and principled approach for obtaining models that can in principle model any combination of network features. In this approach a certain collection of graph features $\{t_i\}$ are constrained to their observed values and the graph is  otherwise assumed to be maximally random where the randomness of the model is measured in terms its Shannon entropy \cite{Cover2012ElementsTheory}. In reference to equilibrium statistical mechanics we will refer such models as canonical and microcanonical ensembles \cite{Bianconi2009EntropyEnsembles} for the cases where constraints are satisfied in expectation and exactly, respectively.  In canonical ensembles the distribution over graphs is given by an exponential of the constrained quantities and therefore such models are also known as exponential random graph models (ERGMs). Whereas in the microcanonical case the ensemble that maximizes the entropy is the ensemble where every configuration that satisfies the given constraints has equal probability. Although maximum entropy models seem to offer an elegant method for constructing random graph models with any desired set of features they are known to be notoriously hard to approach analytically when higher order subgraphs are included in their formulation. Despite these difficulties most efforts for modelling non-trivial subgraph structures have been focused on ERGMs that have counts of triangles and other subgraphs as their parameters \cite{2012ExponentialNetworks,Fienberg2010IntroductionData}. Because such ERGMs are not analytically tractable except in a few isolated instances \cite{Park2004StatisticalNetworks} their study has been mostly been confined to Monte Carlo approaches which themselves suffer from issues of degeneracy and inconsistency \cite{Chatterjee2013EstimatingModels}. 

 In this article we follow the common conception that network motifs are basic building blocks of networks and develop a class of maximum entropy models that is based on constraining counts and distributions of atomic subgraphs \textit{used to construct} the network rather than the counts of subgraphs in the final network. The resulting models can generate networks with a large variety of local structures while remaining analytically solvable for many of their properties. The models we consider fall in the same category as some more recent models that use explicit copies of higher order atomic subgraphs  \cite{Bollobas2011SparseClustering,Newman2009RandomClustering,Miller2009PercolationNetworks,Karrer2010RandomSubgraphs}.


The assumption that networks are formed by atomic subgraphs naturally leads us to consider objects we call subgraph configurations which correspond to the set of atomic subgraphs added to the graph during the generation process. Subgraph configurations are a generalization of hypergraphs where groups of vertices are connected by hyperedges that are not necessarily cliques and can have arbitrary topology.  Models that consider higher order interactions \cite{Battiston2020NetworksDynamics} in the form of cliques have been widely studied before in the form of random hypergraphs \cite{Ghoshal2009RandomApplications}, bipartite models \cite{Guillaume2004BipartiteNetworks,Newman2001RandomApplications} and simplicial complexes \cite{Courtney2016GeneralizedComplexes}. However, the assumption that higher order interactions are cliques is rather restrictive and does not generalize well to directed and/or signed networks. More specifically we focus on maximum entropy ensembles of subgraph configurations given various types of constraints on the counts and the distributions of atomic subgraphs. Starting with the general case we derive expressions for the entropy of canonical and microcanonical ensembles. We also consider a systematic procedure for relaxing constraints placed on the distribution of atomic subgraphs which results in more coarse grained models of varying parametric complexity.

The article is organized as follows. In Sec.\ref{SGMCs} we introduce subgraph configurations and related concepts. We then consider canonical (Sec.\ref{canonical}) and microcanonical (Sec.\ref{Mcanonical})  ensembles of subgraph configurations and present general expressions for their entropy. We also discuss several special cases starting with random graph models for graphs with non-trivial local structures (Sec.\ref{rgmotifs}). We then consider models with labelled atoms and their relation to block models (Sec.\ref{SBM}) and multilayer networks (Sec.\ref{Mlayer}). We conclude with a summary of our main results and potential directions of future studies in Sec.\ref{conc}. 

\section{Subgraph configuration models}\label{SGMCs}

\subsection{Isomorphisms, motifs and orbits}
\begin{figure}
    \centering
    \includegraphics[width=0.15\textwidth]{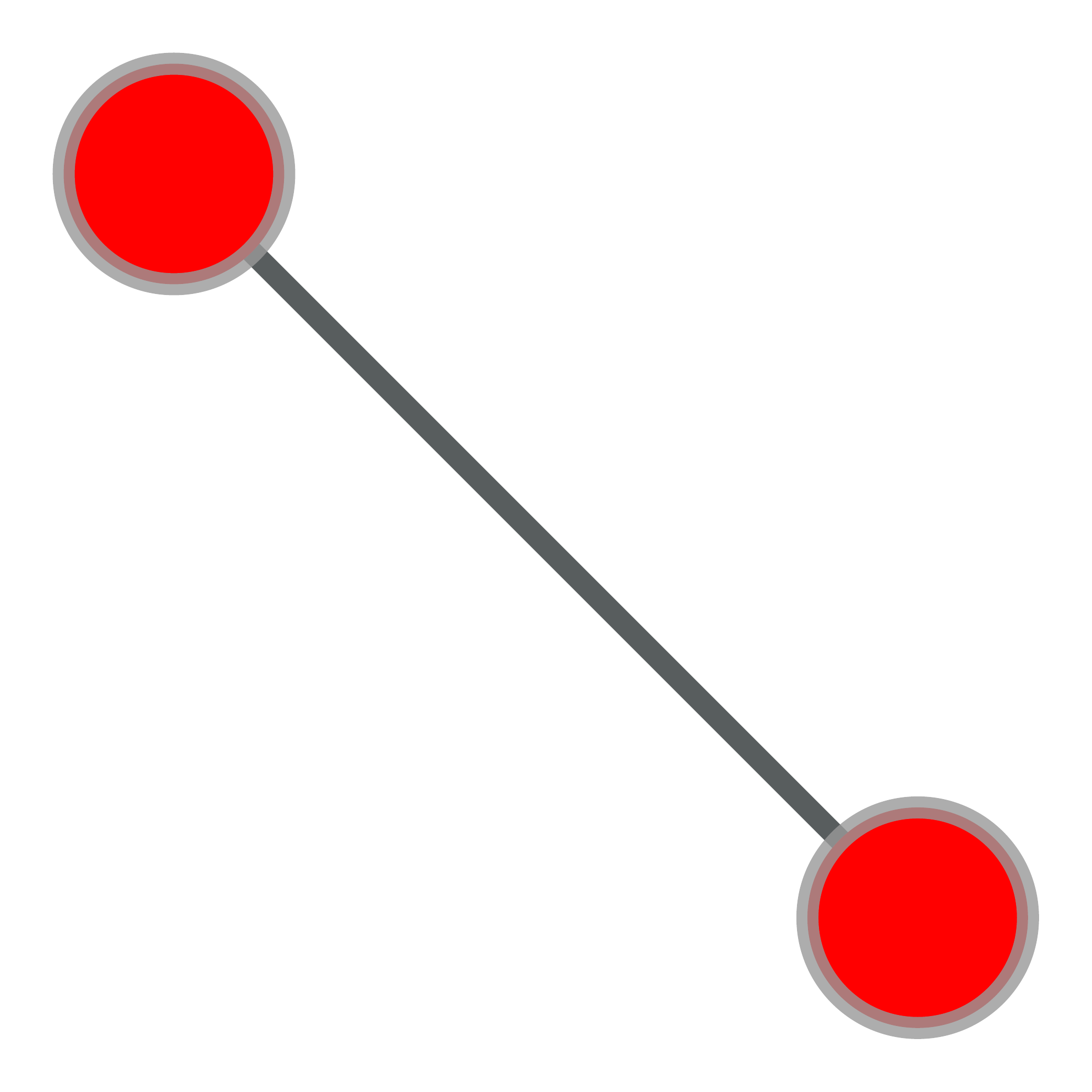}
    \includegraphics[width=0.15\textwidth]{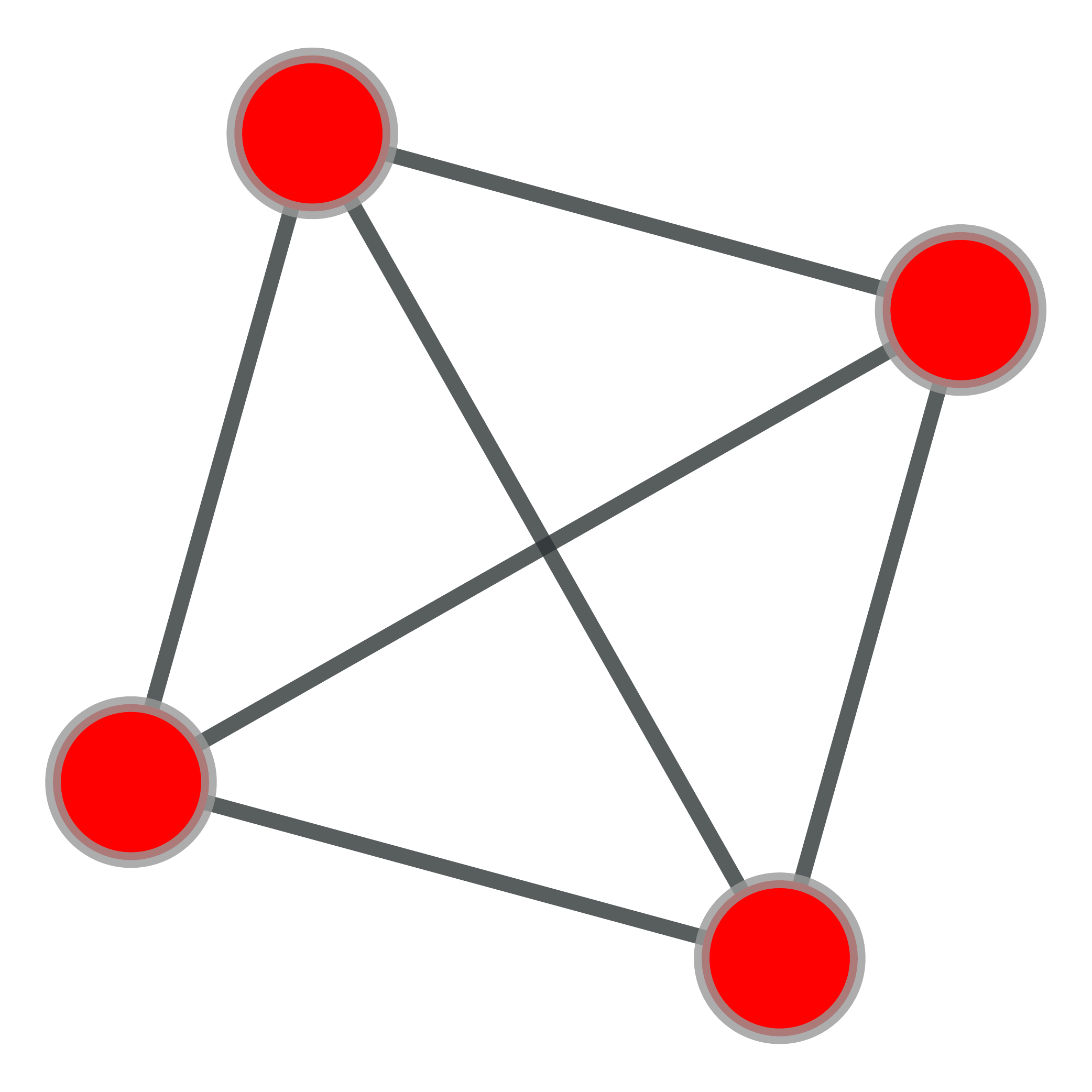}
    \includegraphics[width=0.15\textwidth]{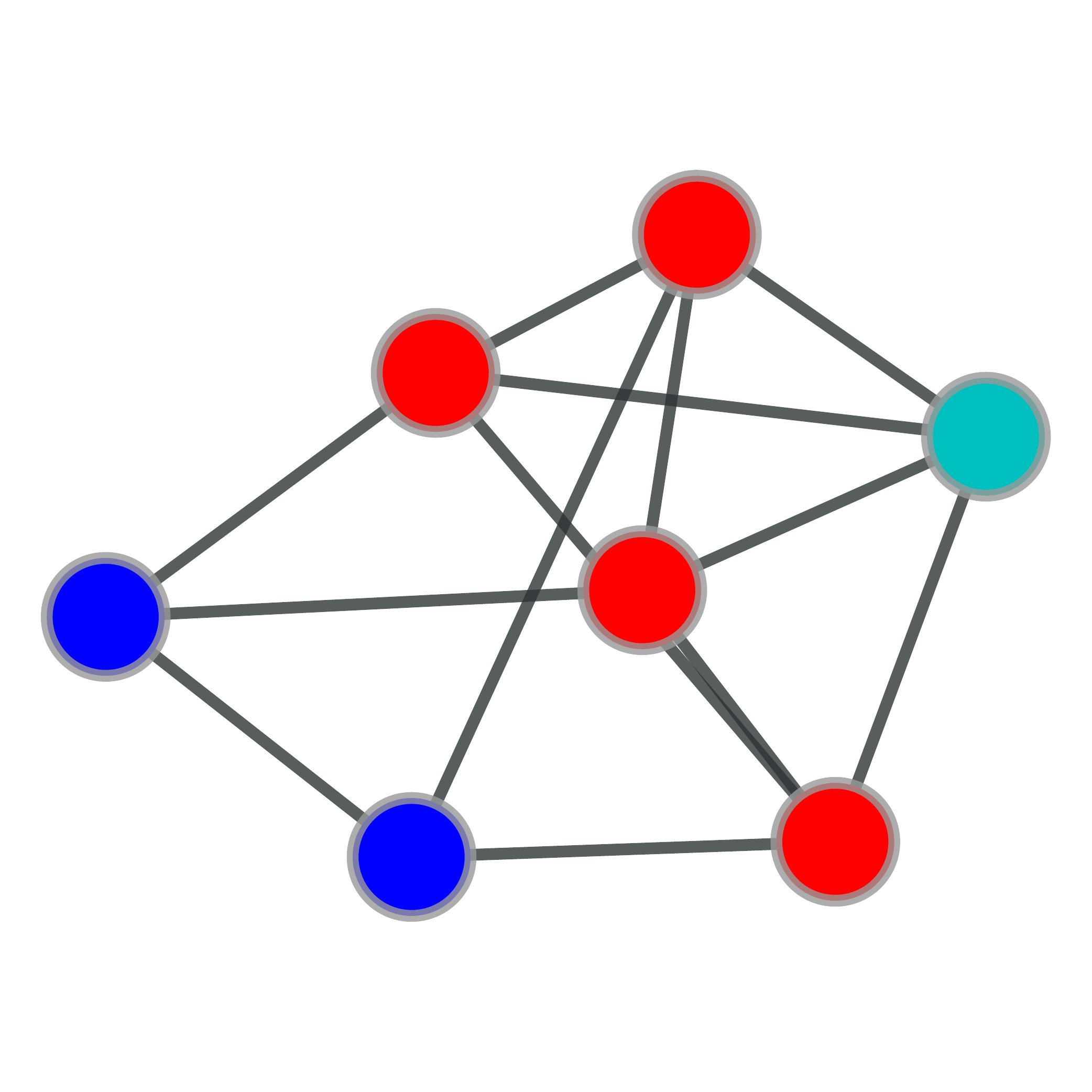}
    \includegraphics[width=0.15\textwidth]{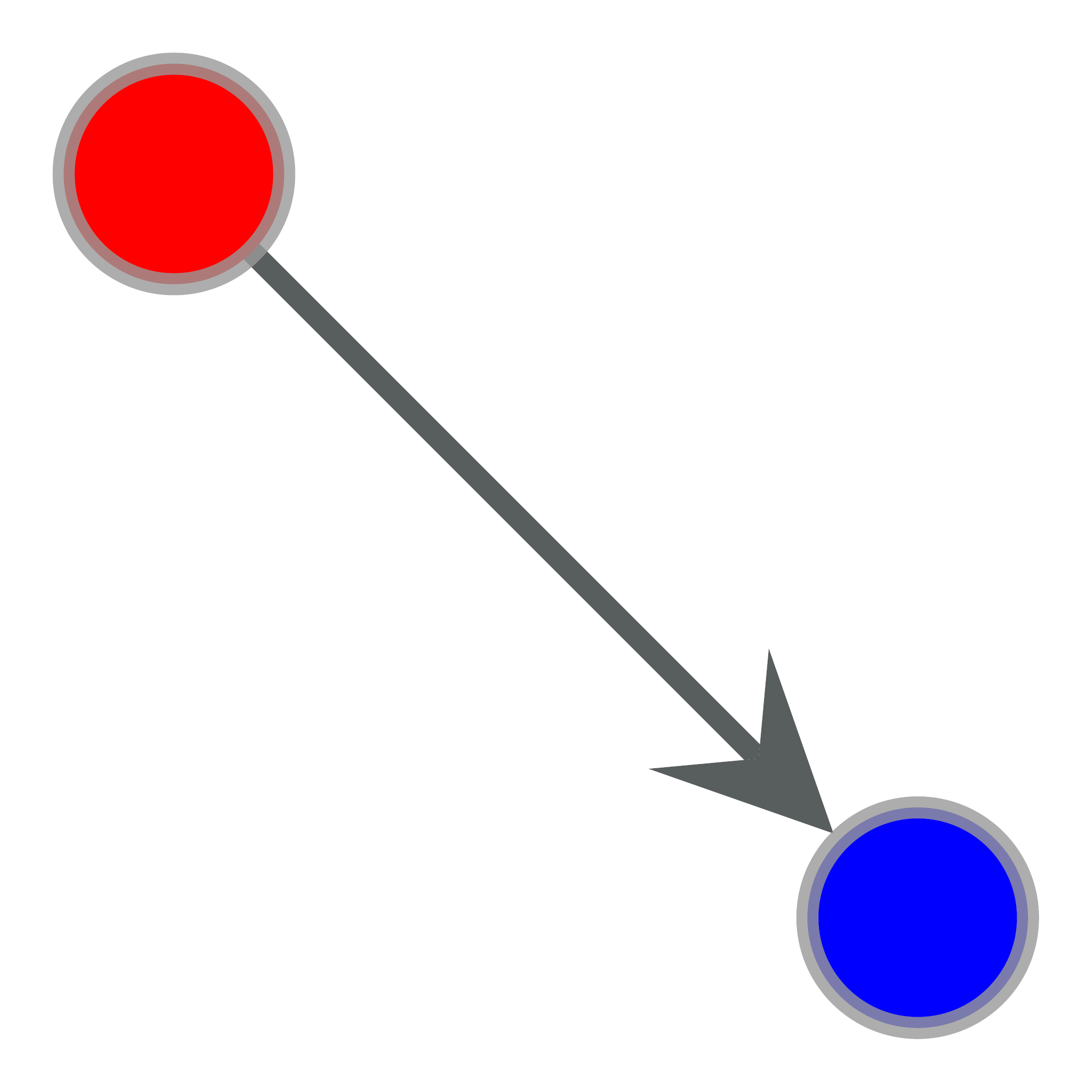}
    \includegraphics[width=0.15\textwidth]{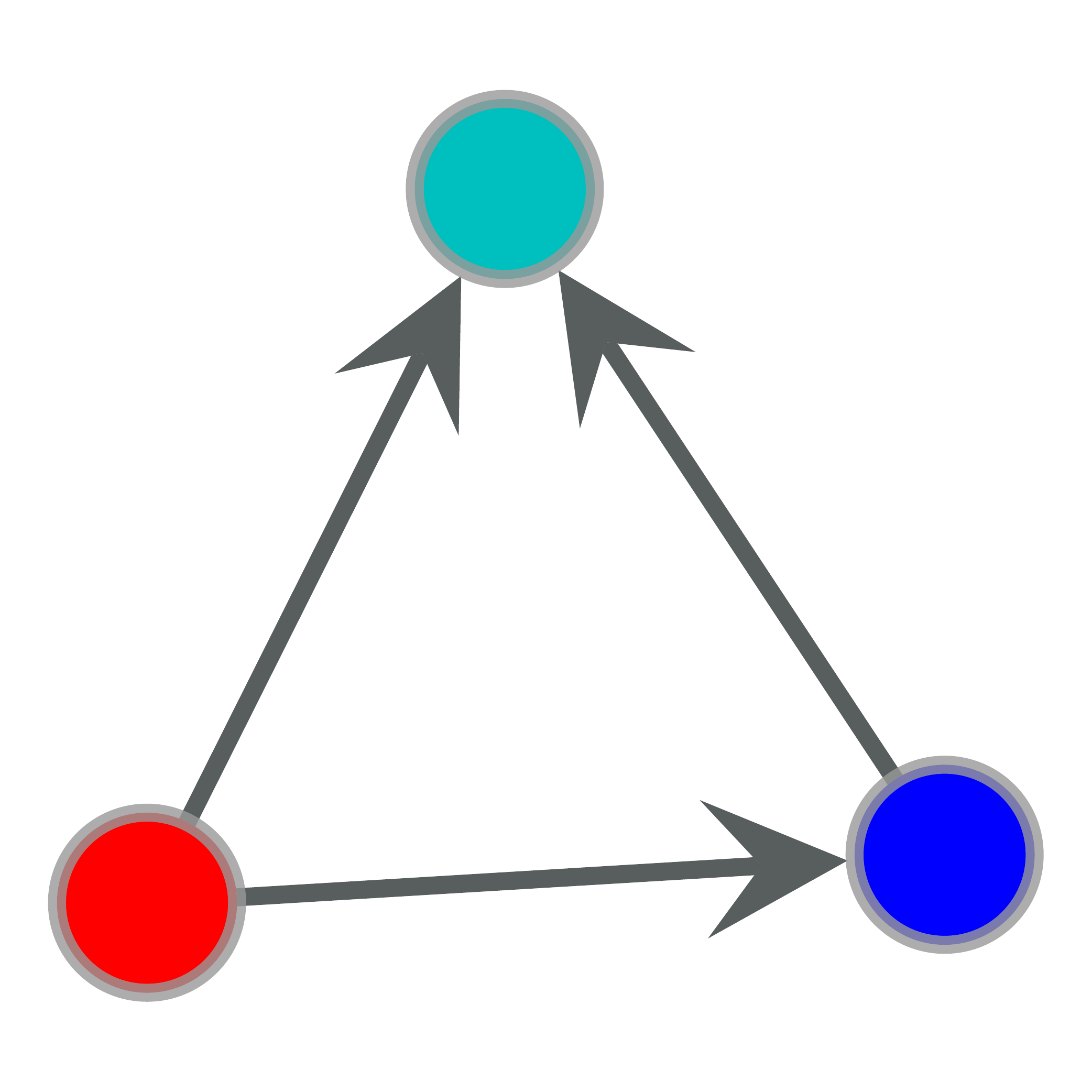}
    \includegraphics[width=0.15\textwidth]{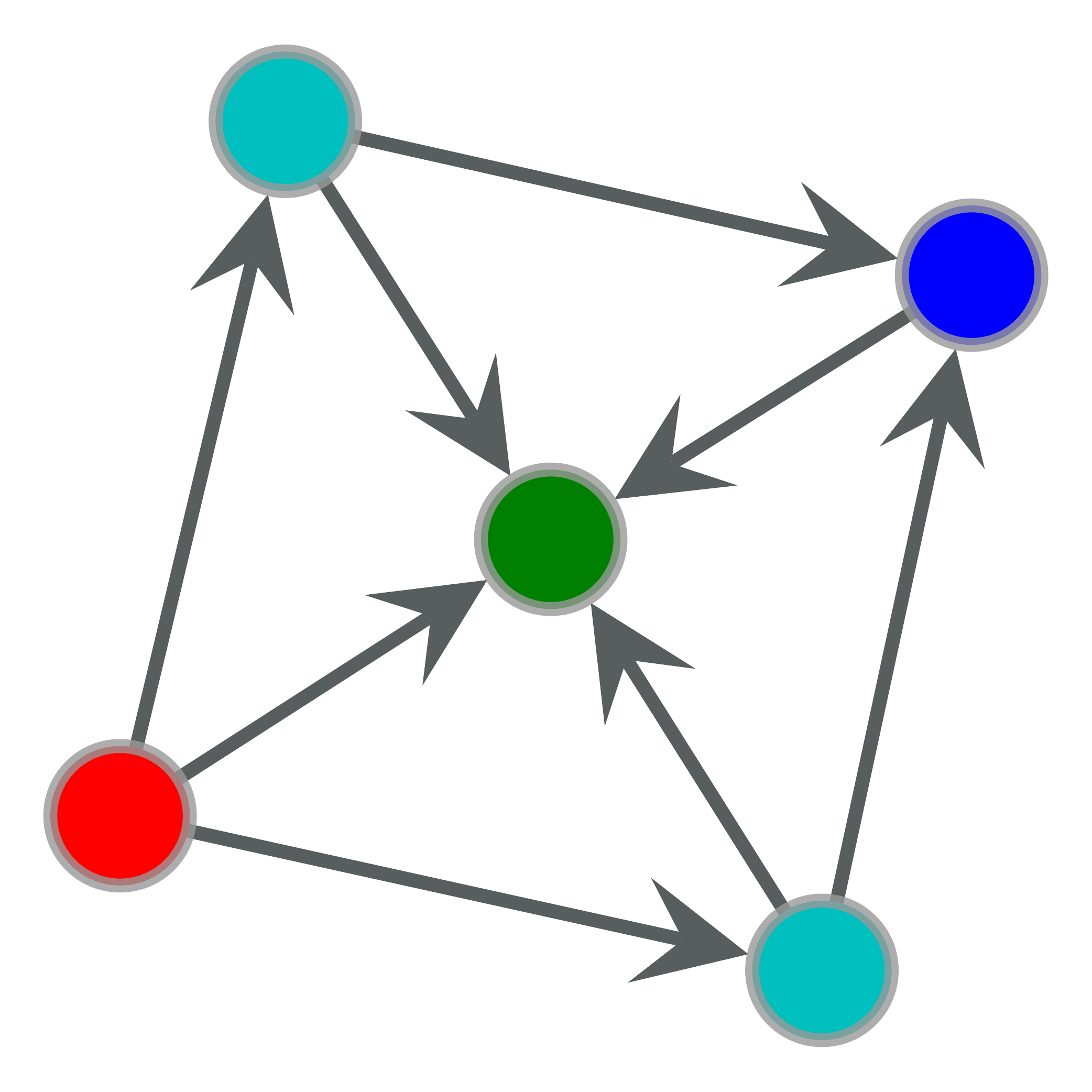}
    \caption{Examples of undirected and directed atoms. Vertex colours indicate the orbits of the atoms.}
    \label{fig:Atoms}
\end{figure}

Before introducing the models we briefly overview key graph theoretical concepts and definitions that are used throughout the text. For a graph $G(E,V)$ we denote its vertex set as $V(G)$ and its edge set as $E(G)$. We sometimes denote the number of vertices/order of $G$ as $|G|$. Symmetry plays an essential role in describing configurations of atomic subgraphs. Two graphs $G$ and $H$ are said to be isomorphic if there exists a bijection $\phi: V(G) \rightarrow V(H)$ such that $(v,v')\in E(G) \iff (\phi(v),\phi(v')) \in E(H)$. If the graphs are directed and/or have labelled (coloured) edges and/or vertices an isomorphism $\phi$ has to also preserve edge directions and labels. If $G=H$ in the above definition $\phi$ is called an automorphism. The automorphisms of $G$ form a group under composition which we denote as $\mathrm{Aut(G)}$. We call the orbits formed by the action of $Aut(G)$ on $V(G)$ the orbits of $G$ and denote the $i^{th}$ orbit of $G$ as $O_{G,i}$. The orbits of a graph are classes of vertices which can be mapped onto each other by vertex permutations that leave the structure of the graph unchanged (See Fig.\ref{fig:Atoms}). Being isomorphic is an equivalence relation of which the equivalence classes we refer to as motifs. We denote motifs using lower-case letters. The automorphism group and orbits of a graph are uniquely determined by its isomorphism class. 

A graph $H$ is said to be a subgraph of $G$ iff $V(H)\subseteq V(G)$ and $E(H)\subseteq E(G)$. Similarly, a $m$-subgraph of $G$ is a subgraph of $G$ that is in the automorphism class $m$. Two subgraphs are said to be distinct unless $E(G)=E(H)$ and $V(G)=V(H)$.


\subsection{Subgraph configurations}
\begin{figure}
    \centering
    \includegraphics[width=0.45\textwidth]{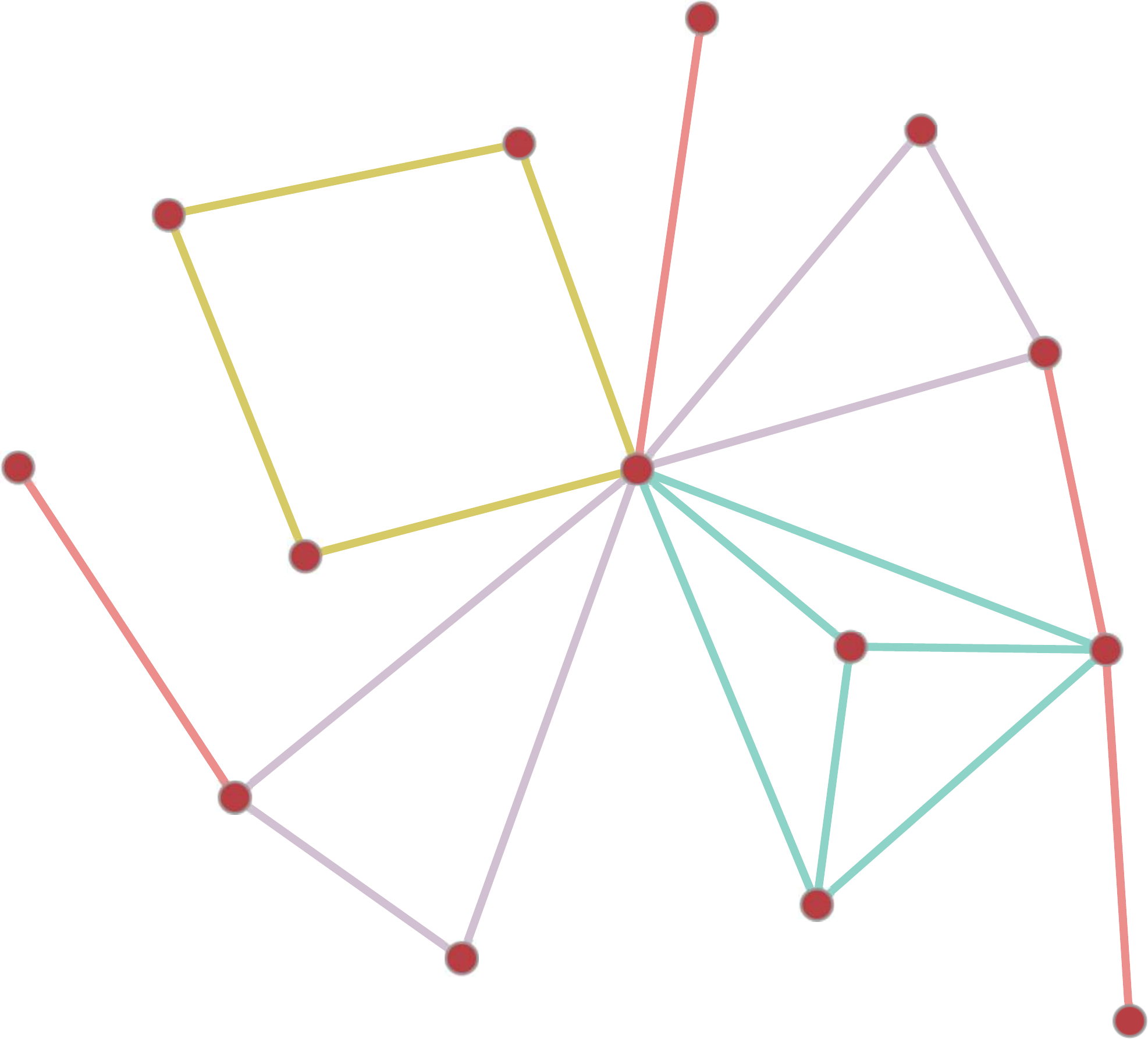}
    \caption{A subgraph configuration consisting of single edges, triangles, a 4-cycle and a 4-clique.}
    \label{fig:Conf}
\end{figure}

A subgraph configuration $C$ on a set of vertices $V$ is a set of subgraphs of the maximally connected graph $K_V$ on $V$. In other words $K_V$ is the graph that contains all possible edges. The specifics of $K_V$ depend on the type of graph under consideration for instance whether it is directed, contains self loops, has multiple layers, is bipartite etc. For an example of a subgraph configuration see Fig.\ref{fig:Conf}.

The set of all $m$-subgraphs of $K_V$ is denoted as $\mathcal{H}_{V,m}$. For undirected simple graphs it follows from the definition of the automorphism group that for each subset of $|m|$ vertices there are $\frac{|m|!}{|\mathrm{Aut(m)}|}$ possible $m$-subgraphs and therefore for a set of $N$ vertices we have:
\begin{equation}
|\mathcal{H}_{N,m}|=\frac{N!}{(N-|m|)!|\mathrm{Aut(m)}|}.
\end{equation}

The set of atoms of a configuration $M(C)$ is the set of motifs occurring in $C$.  In general we will assume $m\in M$ are connected and do not contain multi edges. Similarly, given a set of atoms $M$ an $M$-configuration is one for which $M(C)\subseteq M$. The set of all $M$-configurations on $V$ can be defined as  $\mathcal{C}_{M,V}=\bigtimes_{m\in M}\{0,1\}^{|\mathcal{H}_{V,m}|}$.  We denote the number of $m$-subgraphs in a configuration $C$ as $n_m(C)$.

Subgraph configurations can be described in terms of subgraph tensors that are similar to the adjacency matrices of graphs. Given a configuration $C$ and subgraph $s$ the subgraph tensor $\sigma_s(C)$ is defined as:
\begin{equation}
\sigma_s(C)=
\begin{cases}
      1, & \text{if}\ s \in C \\
      0, & \text{otherwise}
    \end{cases}.
\end{equation}

Subgraph tensors can also be indexed in terms tuples of vertices in analogy with the adjacency matrix. For instance, given an atom  $m$ in terms of a labelled representative the with vertex set $\{1,...,|m|\}$ the subgraph tensor can be defined to be $\sigma(m)_{(v_1v_2...v_m)}(C)=1$ whenever the map $\phi(i)=v_i$ is an isomorphism to some $s\in C$ and 0 otherwise. This implies that $\sigma(m)_{v_1v_2...v_m}(C)=\sigma(m)_{\beta(v_1v_2...v_m)}(C)$ for any permutation $\beta$ that is in $\mathrm{Aut(m)}$. In other words subgraph tensors corresponding to the atom $m$ have to be invariant under 
$\mathrm{Aut(m)}$. This is similar to the condition $A_{ij}=A_{ji}$ for the adjacency matrix of undirected graphs. 

The number of $m$-subgraph in a configuration $C$ can be written in terms of subgraph tensors:
\begin{equation}\label{nm}
    n_m(C)=\sum_{s\in \mathcal{H}_{V,m}}\sigma_s(C).
\end{equation}

For an $M$-configuration one can define its orbit degree sequence $\mathbf{d}_{m,i}(C)(v)$ as the number of $m$-subgraphs in $C$ for which $v$ is in orbit $O_{m,i}$. In terms of $\sigma_S$ $\mathbf{d}_{m,i}(C)(v)$ can be expressed as:

\begin{equation}\label{dmi}
\mathbf{d}_{m,i}(C)(v)=\sum_{s\in \mathcal{H}_{V,m}|v\in O_{m,i}(s)} \sigma_s(C)
\end{equation}
 
The notion of graphicality of (edge) degree sequences also extends to orbit degree sequences.  For instance, a graphical orbit degree sequence has to contain orbits in the right proportions to be graphical: 
\begin{equation}
    \frac{\sum_v d_{m,i}(v)}{|O_{m,i}|}=n_m \forall m,i.
\end{equation}
Throughout this article we shall assume that all orbit degree sequences under consideration are graphical.

\subsection{Subgraph configurations and graphs}
A subgraph configuration $C$ on vertex set $V$ can be projected onto a graph $G=\mathcal{G}(C)$ on $V$ by taking the union of the edges of the subgraphs in $C$. In general the exact form of the projection will depend on the type of graph under consideration i.e. whether it is directed, has multiple layers, admits parallel edges etc. In general we will assume that graphs are simple and hence that configurations are mapped onto graphs by taking the union edge set of the subgraphs in $C$ i.e. $E(G)=\bigcup_{s \in C} E(s)$ which is equivalent to replacing any edges that occur in multiple times in the configuration by single edges in the graph. For the sparse models we study the expected number of such parallel edges is in general $O(1)$. 

We say that a subgraph configuration $C$ covers $G$ if the projection of $C$ is equal to $G$. We denote the set of all $M$-configurations that cover $G$ (or simply $M$-covers of $G$) as $\mathcal{C}_{M}(G)$. The covers of a graph $G$ are exact representations of $G$ in the sense that given any of its covers $G$ can be recovered fully from it. Indeed many widely used graph representations are special cases of subgraph configurations that are covers. For instance the edge list is equivalent to the configuration consisting of all single edges and the adjacency list equivalent to the configuration that contains for each vertex $v$ the star shaped subgraph where $v$ is connected to all its neighbours. As a result given a cover $C$ of $G$ many of its properties can be derived directly from $C$. For instance connected components of the graph coincide with the connected components of its covers. Similarly, every subgraph $s\in C$ also occurs in $G(C)$ hence the projected graph contains at least as many $m$-subgraphs as the $C$ does. This simple fact allows the subgraph structure of graphs to be controlled by changing the counts and types of atoms in subgraph configuration models.  
 
\subsection{Subgraph configuration models and random graphs}

Given a set of atoms $M=\{m\}$ a subgraph configuration model is simply a probability distribution $\mathcal{P}_M(C)$ over the space of subgraph configurations $\mathcal{C}_{M,V}$. While subgraph configuration models can be used to directly model data sets that include higher order interactions such as hypergraphs, directed hypergraphs and simplicial complexes many data sets cover only pairwise interactions i.e. are given in the form of a graph. In such cases a distribution over the space of graphs $\mathcal{G}_V$ can be obtained using the projection defined in the section above. 

Consequently, the distribution over graphs $\mathcal{P}_M(G)$ induced by a subgraph configuration model $\mathcal{P}_M(C)$ is given by: 
\begin{equation}\label{pc}
\begin{split}
P_M(G)&= \sum\limits_{C \in \mathcal{C}_{M,V}} P_M(C)\delta{(G,\mathcal{G}(C))}\\
&=\sum\limits_{C \in \mathcal{C}_M(G)} P_M(C),
\end{split}
\end{equation}
where $\delta(G,\mathcal{G}(C))$ is one whenever $G=\mathcal{G}(C)$ and zero otherwise. In other words the probability of $G$ in the model is given by the total probability of all configurations of which the projection is $G$. 

Subgraph configuration models differ from most other latent state models in that each latent state projects to a single graph. This allows many properties of the model at the graph level to be calculated at the level of configurations.  Another consequence of Eq. \eqref{pc} is that the entropy of a subgraph configuration model $\mathcal{P}_M(C)$ is an upper bound for the entropy of $\mathcal{P}_M(G)$ the distribution it induces on graphs.  
 
\subsubsection{Multi-occupancy subgraph configurations}

It is possible to consider versions of the subgraph configuration models where a configuration can contain multiple copies of the same subgraph. This modification is straightforward and much of the results for the single occupancy and multi-occupancy variants coincide in the sparse setting where the expected number of multiple 'parallel' subgraphs is $o(1)$ for atoms of order higher than 2 and is $O(1)$ for 2 vertex atoms/edges. Hence any modifications to the expressions obtained in this  article in the case of multi-occupancy configurations are dominated by the contribution of two vertex atoms. The case of multi-graph ensembles has been studied extensively before e.g. in Ref.\cite{Peixoto2012EntropyEnsembles}.

\section{Canonical ensembles of subgraph configurations}\label{canonical}
Canonical subgraph configuration ensembles are maximum entropy distributions under constraints given in the form of expectations. In our case the (Shannon) entropy \cite{Cover2012ElementsTheory} of a subgraph configuration model $\mathcal{P}_M(C)$ is defined as:
\begin{equation}
    S(\mathcal{P}_M(C))=-\sum_{C\in\mathcal{C}_{M,V}}P_M(C)\ln(P_M(C)).
\end{equation}

Given a set of atoms $M$ and a set constraints on the expectations of a given set of features $\{t_1,t_2...t_n\}$ the maximum entropy distribution over $\mathcal{C}_{M,V}$ takes the well known exponential form: 

\begin{equation}\label{exp}
    P_M(C)=\frac{1}{Z}\exp(-\sum_i \alpha_i t_i(C)),
\end{equation}
where $Z=\sum_{C  \in \mathcal{C}_{M,V}}\exp(-\sum_i \alpha_i t_i(C))$ is a normalizing constant known as the partition function. As a result canonical ensembles of subgraph configurations are generalizations of ERGMs to hypergraphs where the topology of admissible hyper-edges is given by $M$. 

In general the enumeration of the partition function is a major technical challenge in obtaining analytical results for ERGMs that include higher order interactions in the form of subgraph counts. In our case though we restrict ourselves to features that can be expressed as linear combinations of subgraph tensors resulting in models that are analytically tractable.

\subsection{Canonical ensembles with given expected atomic subgraph counts}

Given a set of atoms $M=\{m\}$ the simplest type of constraint that can be placed on a canonical ensemble is to fix the expected counts of atoms $n_m$ for $m \in M$:

\begin{equation}
E(n_m)=c_m \text{ for } m\in M.
\end{equation}

In general we will focus on sparse graphs and hence assume that $c_m=O(N)$. In such models atoms are distributed uniformly over the vertices of the network an reduce to the Erd\"os-Reny\'i random graph $G(N,p)$ when $M$ only contains the single edge atom.  
Combining Eq. \eqref{nm} and \ref{exp} we obtain that each $m$-subgraph in $\mathcal{H}_{N,m}$ occurs independently with probability:
\begin{equation}
    p_m=\frac{e^{-\lambda_m}}{1+e^{-\lambda_m}}.
\end{equation}

Imposing the constraints on the $\mathbf{n}_m$ we have:

\begin{equation}
    \mathbf{p}_m |\mathcal{H}_{V,m}| = \mathbf{c}_m.
\end{equation}

Consequently, the entropy can be written as sum over sum over subgraphs:
\begin{eqnarray}
    S(M,\mathbf{c}_m)&=& \sum_{m\in M}|\mathcal{H}_{
    V,m}| h\Big(\frac{c_m}{|\mathcal{H}_{V,m}|}\Big)\\
    &\simeq& \sum_{m\in M}\big(c_m-c_m\ln\Big(\frac{c_m}{|\mathcal{H}_{V,m}|}\Big)\big),
\end{eqnarray}

where $h(p)=-p\ln(p) - (1-p)\ln(1-p)$ is the binary entropy. 

\subsection{Canonical ensembles with given expected atomic degree sequence}\label{dccm}

Given a set of atoms $M=\{m\}$ the constraints for the atomic degree sequence can be written as:

\begin{equation}
    E(d_{m,i}(v))=k_{m,i}(v),
\end{equation}
for all $m$ in $M$ and their orbits $i$.

Combining Eq.\ref{dmi} and Eq.\ref{exp} results in a expression that can be factorized over the contributions of individual subgraphs:
\begin{equation}
\begin{split}
P(C)&=\frac{1}{Z}\exp{\big(-\sum_{m,i}\sum_{v}\lambda_{m,i}(v)d_{m,i}(C)(v)\big)}\\
&=\frac{1}{Z}\exp{\big(-\sum_{m,i}\sum_{v}\sum_{s\in\mathcal{H}_{V,m}|v\in O_{m,i}(s)}\lambda_{m,i}(v)\sigma_s(C)\big)}\\
&=\frac{1}{Z}\prod_{s\in\mathcal{H}_{V,m}}
\exp{\big(-\sum_{v\in O_{m,i}(s)}\lambda_{m,i}(v)\sigma_s(C)\big)}.
\end{split}
\end{equation}
Where the partition function is given by:
\begin{equation}
Z=\prod_{m\in M}\prod_{s\in\mathcal{H}_{V,m}} ({1+e^{-\sum_{v\in O_{m,i}(s)}\lambda_{m,i}(v)}}).
\end{equation}

Hence in this ensemble every $m$-subgraph $s\in\mathcal{H}_{V,m}$ occurs independently with probability:

\begin{equation}\label{pdmi}
    p_s= \frac{e^{-\sum_{v\in O_{m,i}(s)}\lambda_{m,i}(v)}}{1+e^{-\sum_{v\in O_{m,i}(s)}\lambda_{m,i}(v)}}.
\end{equation}
As a result the entropy of the ensemble in terms of the binary entropy $h(p_s)$:
\begin{equation}
    S(M,\mathbf{n}_m)=\sum_{m \in M} \sum_{s\in\mathcal{H}_{N,m}} h(p_s).
\end{equation}
Similarly, the expectations of atomic degrees can be written as:
\begin{equation}\label{Edmi}
    E(d_{m,i}(v))=\sum_{s\in \mathcal{H}_{N,m}|v\in O_{m,i}(s)} p_s.
\end{equation}

\subsubsection{The sparse limit}
Unfortunately, Eq.\ref{pdmi} and \ref{Edmi} generally do not have a closed form solution. However, if we assume that the $p_s\ll1$ for all subgraphs $s$ and $N\gg1$ we have: $p_s\simeq e^{-\sum_{v\in s}\lambda_{m,i}(v)}$. For $N\gg1$ the expected counts $\Bar{n}_m$ can be approximated as:
\begin{equation}
\begin{split}
    \bar{n}_m&=\sum_{s\in\mathcal{H}_{V,m}} e^{-\sum_{v\in O_{m,i}(s)}\lambda_{m,i}(v)}\\
    &=\frac{\prod_i|O_{m,i}|!}{|\mathrm{Aut(m)}|}\sum_{t\in {V \choose |m|}} \sum_{o \in {t \choose \mathbf{|O_{m,i}|}}}e^{-\sum_{v\in t}\lambda_{o(v)}(v)}\\
    &= \frac{1}{|\mathrm{Aut(m)}|}\prod_i \Big(\sum_{v}e^{-\lambda_{m,i}(v)}\Big)^{|O_{m,i}|} (1+O(1/N)),
\end{split}
\end{equation}
where in the first step the sum over $m$-subgraphs is converted to a sum over $|m|$-subsets of vertices $(t)$ and grouping these subgraphs according the orbit assignments ($o$) of the vertices in $t$. For each such orbit placement $o$ there are $\frac{\prod_i|O_{m,i}|!}{|\mathrm{Aut(m)}|}$ $m$-subgraphs compatible with $o$. The final expression is obtained by converting this to a sum over $|m|$-tuples of vertices which for $|m|\ll N$ can approximated by the sum over $V^{|m|}$.

Substituting this into Eq. \eqref{Edmi} one gets:
\begin{equation}
\begin{split}
    k_{m,i}(v)&= e^{-\lambda_{m,i}(v)}\sum_{s|v\in O_{m,i}(s)}\prod_j \prod_{v'\in O_{m,j}(s)|v'\neq v} e^{-\lambda_{m,j}(v')}\\
    &\simeq e^{-\lambda_{m,i}(v)} \frac{|O_{m,i}|\Bar{n}_m}{ \sum_{v'} e^{-\lambda_{m,i}(v')}}.
\end{split}    
\end{equation}

Solving the system of equations we obtain:

\begin{equation}\label{canonicalP}
    p_s=\bar{n}_m |\mathrm{Aut(m)}|\prod_{i,v|v\in O_{m,i}(s)}\frac{k_{m,i}(v)}{|O_{m,i}|\bar{n}_m}, 
\end{equation}
where $\bar{n}_m$ is the expected number of total $m$-subgraphs. 
In order for the model to be well defined we need $p_s\leq 1$ which implies that:
\begin{equation}\label{cond}
\prod_{i,v|v\in O_{m,i}(s)} k_{m,i}(v)\leq \bar{n}_m^{|m|-1}\prod_i|O_{i,m}|^{|O_{i,m}|}/|\mathrm{Aut(m)}|,
\end{equation}
for all $s$ in $\mathcal{H}_{V,m}$. This can also be expressed in terms of the average degree $\langle {k}_{m,i}\rangle=N^{-1}\sum_v k_{m,i}(v)$ using the identity: $|O_{m,i}|\bar{n}_m=\sum_v k_{m,i}=N\langle {k}_{m,i}\rangle$.  Although because of the potential presence of multiple orbits there are many ways to satisfy this constraint, the simplest constraint is to set: $max_{\{k_{m,i}(v)\}}<|O_{m,i}|\bar{n}_m^{(|m|-1)/|m|}|\mathrm{Aut(m)}|^{1/|m|}$. Note that this reduces to the familiar condition $k_{max}<\sqrt{\langle k\rangle N}$ when $m$ is the single edge. Because in most practical cases the model is likely to contain atoms consisting of edges we shall assume that $max\{k_{m,i}(v)\} \ll \sqrt{\langle k_{m,i}\rangle N}$ which guarantees Eq. \eqref{cond}.

Using the expansion $h(x)=-x\ln(x) +x -\sum_{l=1}^\infty \frac{x^{l+1}}{l(l+1)}$ for the binary entropy we obtain the following expression for the entropy:
\begin{equation}\label{Centropy}
\begin{split}
    &S(M,\mathbf{d}_{m,i})=\sum_m\Bigg[\bar{n}_m-\bar{n}_m\ln(|\mathrm{Aut(m)}|\bar{n}_m)\\&+\sum_i \bar{n}_m|O_{m,i}|\ln(\bar{n}_m|O_{m,i}|) -\sum_i\sum_v k_{m,i}(v) \ln(k_{m,i}(v)) \\&-\sum_{l=1}^\infty\frac{|\mathrm{Aut(m)}|^{l}}{l(l+1)} \frac{\bar{n}_m^{l+1}}{\Big(\prod_j(|O_{m,j}|\bar{n}_m)^{|O_{m,j}|}\Big)^l }\prod_i\Big(\frac{\langle k_{m,i}^{l+1}\rangle}{\langle k_{m,i}\rangle}\Big)^{|O_{m,i}|}\Bigg].
    \end{split}
\end{equation}
Eq. \eqref{Centropy} generalizes many known formulas for degree corrected canonical ensembles including directed and undirected graphs \cite{Bianconi2009EntropyEnsembles,Bender1978AsymptoticFunctions} to which it reduces when $M$ contains only the single edge atom. Note that in our formulation the difference between undirected and directed graphs emerges simply through the difference of the automorphism groups of undirected and directed single edge atoms.

\subsubsection{Orbit aggregation}
Degree corrected subgraph configuration models require a sequence of length $N$ to be specified for each orbit of the atoms in $M$. Hence, for large $M$ the number of parameters can quickly become excessive. Therefore we now present models where constraints on the atomic degrees are relaxed by combining the degrees of two or more orbits. Aggregating orbits results in models with lower parametric complexity which is especially relevant in the context of statistical inference where in general the goal is to obtain a model that fits the data well without requiring an excessive number of parameters. 

Aggregating two different orbits and specifying their total degree results in a constraint in the form:
\begin{equation}\label{aggconst}
    E(d_{m_1,i}(v)+d_{m_2,j}(v))=k(v).
\end{equation}
Note that $m_1$ and $m_2$ in the above expression can be the same atom. We also assume that $\bar{n}_{m_1}$ and $\bar{n}_{m_2}$ are given. The same derivation as in Sec.\ref{dccm} can be carried over to this case and we get the analogue of Eq.\ref{pdmi}:

\begin{equation}\label{psagg}
    p(s)= \frac{e^{-\lambda_m-\sum_{v\in O_{m,i}(s)}\lambda_{m,i}(v)}}{1+e^{-\lambda_m-\sum_{v\in O_{m,i}(s)}\lambda_{m,i}(v)}},
\end{equation}
for $m_1$ and $m_2$-subgraphs where  $\lambda_{m_1,i}(v)=\lambda_{m_2,j}(v)=\lambda(v)$.  In the classical limit imposing the constraints in Eq.\ref{aggconst} we obtain:
\begin{equation}
    k(v)\approx e^{-\lambda(v)} \frac{|O_{m_1,i}|\bar{n}_{m_1}+|O_{m_2,j}|\bar{n}_{m_2}}{ \sum_{v'} e^{-\lambda(v')}}
\end{equation}

Consequently, the probability of a subgraph $s$ of type  $m_1$ or $m_2$ has the same form as Eq. \eqref{canonicalP} with the following substitutions:
\begin{eqnarray}\label{canonAgg}
   k_{m_1,i}(v)&=& \frac{k(v)|O_{m_1,i}|\bar{n}_{m_1}}{|O_{m_1,i}|\bar{n}_{m_1}+|O_{m_2,j}|\bar{n}_{m_2}},\\
   k_{m_2,j}(v)&=& \frac{k(v)|O_{m_2,j}|\bar{n}_{m_2}}{|O_{m_1,i}|\bar{n}_{m_1}+|O_{m_2,j}|\bar{n}_{m_2}}.
\end{eqnarray}
Note that orbit aggregation is commutative and the above expressions generalize to the case where more than two orbits are aggregated as well as the case where multiple (disjoint) sets of orbits are combined. The entropy for the ensemble can be obtained by making the above substitutions in Eq. \eqref{Centropy}.

\section{Microcanonical ensembles}\label{Mcanonical}
In this section we consider subgraph configuration models which obey hard constraints on atomic subgraph counts and degrees. For hard constraints the maximum entropy ensemble is the one where all configurations satisfying the constraints are equiprobable and configurations not satisfying the constraints have zero probability. Consequently, in microcanonical ensembles we have $P(C)=1/\Omega$  for any $C$ satisfying the constraints and the entropy is given by $S=\ln(\Omega)$, where $\Omega$ is the total number of configurations that satisfy the given constraints.

\subsection{Microcanonical ensemble with fixed atomic subgraph counts}
Given a set of atoms $M$ and corresponding subgraph counts $n_m$ there are:
\begin{equation}
\Omega(M,\mathbf{n}_m)=\prod_{m\in M} {|\mathcal{H}_V(m)| \choose n_m},
\end{equation}
such configurations. For instance when $M$ consists of the edge and triangle motifs the microcanonical ensemble with $e$ edges and $t$ triangles induces a probability distribution over all graphs that can be constructed using $e$ edges and $t$ triangles. Note that however the distribution over such graphs (Eq.\eqref{pc}) is not uniform since the probability of a given graph $G$ is proportional to the number of different ways it can be constructed using $e$ edges and $t$ triangles.

\subsection{Microcanonical ensembles with fixed atomic degree distribution}\label{sgcm}

The entropy of micro canonical ensembles can be derived following two approaches. We first consider an analytic approach that is based on imposing hard constraints on the canonical ensemble and later a combinatorial approach that is based on the generative model by Karrer and Newman  \cite{Karrer2010RandomSubgraphs}.  

\subsubsection{Analytic approach}
We first base our treatment of the microcanonical ensembles of subgraph configurations with a given atomic degree sequence on the conjugate canonical ensemble. The entropy of the microcanonical ensemble $S_{mc}$ can be expressed in terms of the entropy of the canonical ensemble $S_c$ as:
\begin{equation}\label{ctomc}
S_{mc}=S_{c}-\Omega,  
\end{equation}
where $\Omega=-\sum_{m\in M}\sum_{i\in O(m)}\sum_{v\in V}\ln(\pi(d_{m,i}(v)))$ and $\pi(k)=\frac{k^k}{k!}e^{-k}$ are Poisson probabilities. A proof of this result for graphs is given in \cite{Bianconi2009EntropyEnsembles} and for simplicial complexes in \cite{Courtney2016GeneralizedComplexes} hence we omit the details of the generalization to subgraph configurations. Applying the above results yields the following expression for the entropy of the microcanonical ensemble:

\begin{equation}\label{Smcstat}
\begin{split}
    &S_{mc}(M,\mathbf{d}_{m,i})=\sum_m\Bigg[n_m-n_m\ln(n_m) -n_m\ln(|\mathrm{Aut(m)}|)\\&+\sum_i n_m|O_{m,i}|(\ln(n_m|O_{m,i}|)-1) -\sum_i\sum_v \ln(d_{m,i}(v)!) \\&-\sum_{l=1}^\infty\frac{|\mathrm{Aut(m)}|^{l}}{l(l+1)} \frac{n_m^{l+1}}{\Big(\prod_j(|O_{m,j}|n_m)^{|O_{m,j}|}\Big)^l}\prod_i\Big(\frac{\langle d_{m,i}^{l+1}\rangle}{\langle  d_{m,i}\rangle}\Big)^{|O_{m,i}|}\Bigg].
\end{split}
\end{equation}

\subsubsection{Combinatorial approach}

Another model that is closely related to subgraph configurations is the model introduced independently by Newman \cite{Newman2009RandomClustering} and Miller \cite{Miller2009PercolationNetworks} which generalizes the classical (edge) configuration model to the case where in addition to the edge degree vertices are also assigned triangle degrees. The model was later generalized by Newman and Karrer \cite{Karrer2010RandomSubgraphs} to allow for arbitrary atomic subgraphs. 

In Ref.\cite{Karrer2010RandomSubgraphs} the authors describe their model in terms of a generating process analogous to the stub matching process for the edge configuration model. In this process given a set of atoms $M$ and a corresponding atomic degree sequence $\mathbf{d}_{m,i}$ one attaches to every vertex atomic stubs reflecting its atomic degrees. Atomic stubs are partial subgraphs such as half edges in the case of edges and corners in the case of triangles. Though in general one might have different kinds stubs corresponding to the orbits of same atom. A network is then generated by matching stubs corresponding to the same atom $m$ in appropriate combinations uniformly at random and connecting them to form an $m$-subgraph until all stubs are exhausted. For instance, if $M$ consists of triangles and edges one matches the pairs of edge stubs and triples of triangle stubs. This process samples all possible matchings uniformly. However, the process allows for stubs attached to the same vertex to be matched to each other resulting in a subgraph that is a vertex contraction of the original atom. For instance, the vertex contraction of an edge creates a self loop and the vertex contraction of a triangle creates two parallel edge with a self loop on one of its vertices. Moreover, the matching process also allows multiple (parallel) copies of the same subgraph to be created. If one excludes these cases by restarting the generating process whenever they occur every subgraph configuration with atomic degree sequence $\mathbf{d}_{m,i}$ is formed with equal probability. 

We proceed with the calculation of $P(C|\mathbf{d}_{m,i})$. Following the construction of Newman and Karrer \cite{Karrer2010RandomSubgraphs} we first consider the number $\Omega(\mathbf{d}_{m,i})$ of possible stub matchings given an atomic degree sequence $\mathbf{d}_{m,i}$. Note that the matching processes for different $m$ are independent. The number of possible matchings for a given atomic degree sequence $\mathbf{d}_{m,i}$ is: 
\begin{equation*}
\Omega(\mathbf{d}_{m,i})=\prod_m\frac{\mu_m^{n_m}}{n_m!}\prod_i\frac{ (|O_{m,i}|n_m)!}{( |O_{m,i}|!)^{n_m}\prod_{v}d_{m,i}(v)!}.
\end{equation*}
Here $|O_{m,i}|n_m!$ is the number of arrangements of stubs of type $m,i$ and the factors $|O_{m,i}|!^{n_m}$ and $\prod_{v}d_{m,i}(v)!$ account for equivalent arrangements of the stubs. Finally, $n_m!$ accounts for the possible rearrangements of the subgraphs and  $\mu_m=\prod_i|O_{i,m}|!/|\mathrm{Aut(m)|}$ is the number of distinct $m$-subgraphs that can be formed given the orbit memberships of its vertices. For instance, there is a only one possible way a triangle can be formed on 3 vertices whereas there are 3 different ways a 4-cycle can be formed on 4 vertices. Note that both motifs have only one orbit. The terms involving $|O_{i,m}|!$ cancel out and one has:
\begin{equation}\label{stubS}
\Omega(\mathbf{d}_{m,i})=\prod_m\frac{1}{|\mathrm{Aut(m)}|^{n_m}n_m!}\prod_i\frac{ (|O_{m,i}|n_m)!}{\prod_{v}d_{m,i}(v)!}.
\end{equation}

However Eq. \eqref{stubS} does include cases where two or more stubs corresponding to the same vertex are matched together. The probability that none of the stubs attached to a given vertex $v$ having atomic degree $d_{m,i}(v)$ are matched together is given by: 
\begin{equation}
\begin{split}
    &P_{c}(d_{m,i}(v))= \prod_m \frac{n_m!}{(n_m-\sum_i d_{m,i}(v))!}\\&\times\prod_i \frac{(|O_{m,i}|n_m-d_{m,i}(v))!|O_{m,i}|^{d_{m,i}(v)}}{(|O_{m,i}|n_m)!}.
\end{split}
\end{equation}
Expanding the above expression using Stirling's approximation we get:

\begin{equation}
\begin{split}
&\ln(P_{c}(d_{m,i}(v)))=-\frac{1}{n_m}{\sum_i d_{m,i} \choose 2}\\&+\sum_i \frac{1}{|O_{m,i}|n_m}{d_{m,i} \choose 2} +O(\frac{1}{N^2}).    
\end{split}
\end{equation}
Assuming independence and summing over all vertices we obtain the following overall correction factor: 
\begin{equation}\label{sloop}
\begin{split}
 \ln(P_{c}(\mathbf{d}_{m,i}))&= -\frac{1}{2}\sum_m\Big[|m|(\frac{\langle (\sum_i d_{m,i})^2\rangle}{\langle \sum_i d_{m,i}\rangle}-1) \\&-\sum_i \big(\frac{\langle d_{m,i}^2\rangle}{\langle d_{m,i}\rangle}-1\big)  \Big],   
\end{split}
\end{equation}
where $\langle d_{m,i}\rangle=\sum_v d_{m,i}(v)/N$.
Although the independence assumption does not hold exactly; the dependence in general is weak and the independence assumption is known to produce results that are consistent with more rigorous analyses for sparse simple graphs i.e. when $M$ consists of only the single edge motif \cite{Bianconi2009EntropyEnsembles}. 

Even after discounting by the above factor we still are left with the possibility that the same subgraph is created multiple times by the matching process. To account for these cases we consider the probability that at least $2$ copies of a certain $m$-subgraph $s$ being created during the matching process. For an $m$-subgraph with orbits $O_{m,i}(s)$ we have:  
\begin{equation}\label{p2}
\begin{split}
    P_2(s)=& \frac{|\mathrm{Aut(m)}|^2}{2}\frac{n_m!} {(n_m-2)!}\prod_j\frac{(|O_{m,j}|(n_m-2))!}{(|O_{m,j}|n_m)!}\\&\times\prod_{i,v|v\in O_{m,i}(s)} \frac{d_{m,i}(v)!} {(d_{m,i}(v)-2)!}. 
\end{split}
\end{equation}
Where we assume that $d_{m,i}(v)\geq2$ for all vertices of $s$, since $P_2(s)=0$ otherwise. The probability of there being no multiple copies of $s$ is simply $1-P_2(s)$. Assuming independence between subgraphs we have:

\begin{equation}\label{pml}
\begin{split}
&\ln(P_{ml}(\mathbf{d}_{m,i}))=\sum_{s\in\mathcal{H}_N(m)}\ln(1-P_2(s))\\
     \simeq&-\frac{|\mathrm{Aut(m)}|n_m^2}{2}\prod_i\frac{1}{(n_m|O_{m,i}|)^{|O_{m,i}|}}\Big(\frac{\langle d_{m,i}^2\rangle}{\langle d_{m,i}\rangle}-1\Big)^{|O_{m,i}|}.
     \end{split}
\end{equation}
Where we assumed $P_2(s)\ll1$ so that $\ln(1-P_2(s))\simeq-P_2(s)$. Note that when $n_m=O(N)$, $\ln(P_{ml})$ scales as $1/N^{|m|-2}$. Hence in general the effect of $P_{ml}$ on the entropy can be neglected for atoms of order higher than 2. 

Finally, $S_{mc}(\mathbf{d}_{m,i})=\ln\big(\Omega(\mathbf{d}_{m,i})P_{c}(\mathbf{d}_{m,i})P_{ml}(\mathbf{d}_{m,i})\big)$ and combining the factors given in Eq.'s \ref{stubS}, \ref{sloop} and \ref{pml} we obtain the following expression for the entropy:
\begin{equation}\label{Cent}
\begin{split}
&S_{mc}(M,\mathbf{d}_{m,i})
=\sum_m\Bigg[-\ln(n_m!)-n_m\ln(|\mathrm{Aut(m)}|)\\&+\sum_i\big[\ln((|O_{m,i}|n_m)!)-\sum_v \ln(d_{m,i}(v)!)\big]\\
&-\frac{|\mathrm{Aut(m)}|n_m^2}{2}\prod_i\frac{1}{(n_m|O_{m,i}|)^{|O_{m,i}|}}\Big(\frac{\langle d_{m,i}^2\rangle}{\langle d_{m,i}\rangle}-1\Big)^{|O_{m,i}|}\\
&-\frac{1}{2}\Big[|m|(\frac{\langle (\sum_i d_{m,i})^2\rangle}{\langle \sum_i d_{m,i}\rangle}-1) -\sum_i \big(\frac{\langle d_{m,i}^2\rangle}{\langle d_{m,i}\rangle}-1\big)  \Big]\Bigg].\end{split}
\end{equation}
The above reduces to known expressions for entropy of microcanonical ensembles of graphs with a given degree distribution \cite{Bianconi2009EntropyEnsembles,Bender1978AsymptoticFunctions} in both the directed and undirected case when $M$ consists only of the single edge atom. Eq. \eqref{Cent} also agrees with Eq. \eqref{Smcstat} up to $O(\ln(N))$ in the sparse limit assuming that the terms for $l>1$ can be neglected in Eq. \eqref{Smcstat}.

\subsection{Orbit aggregation}\label{mcagg}

As in the canonical case microcanonical ensembles can be relaxed by considering constraints in the form of sums of the orbit degrees:
\begin{equation}
   d(v)= d_{m_1,i}(v)+d_{m_2,j}(v),
\end{equation}
where $m_2$ might be the same as $m_1$. Such ensembles can be obtained by treating the stubs of the combined orbits as a single type during the generation process. In this case one has to replace the two factors corresponding to these orbits in Eq. \eqref{stubS} by:
\begin{equation}\label{mcorbagg}
\frac{ (|O_{m_1,i}|n_{m_1}+|O_{m_2,j}|n_{m_2})!}{\prod_{v}d(v)!}.
\end{equation}

The equivalent of $P_c$ (Eq.\ref{sloop}) can be obtained by assuming that the orbit degree $d_{m_1,i}$ follows a binomial distribution in the interval $[0,k(v)]$ with probability $p_{m_1,i}=\frac{|O_{m_1,i}|n_m}{|O_{m_1,i}|n_m+|O_{m_2,j}|n_{m_2}}$. Which yields  $\langle d_{m,i}\rangle=p_{m,i}\langle d\rangle$ and $\langle d_{m,i}^2\rangle=p_{m,i}^2\langle d(d-1)\rangle+p_{m,i}\langle d\rangle$. Similarly the correction factor for multiple subgraphs (Eq. \eqref{pml}) can be obtained by replacing the factor corresponding to $O_{m_1,i}$ by:
\begin{equation}
    \frac{1}{(|O_{m_1,i}|n_{m_1}+|O_{m_2,j}|n_{m_2})^{|O_{m_1,i}|}}\Big(\frac{\langle d^2\rangle}{\langle d\rangle}-1\Big)^{|O_{m_1,i}|}.
\end{equation}
The final expression for the entropy can be obtained by making the same substitutions for $(m_2,j)$ as well. 

As in the canonical ensemble the above procedure generalizes to the case where one combines multiple (disjoint) sets of orbits. The equivalent of Eq. \eqref{Smcstat} be obtained using the substitutions given in Eq. \eqref{canonAgg} and using the combined degree in Eq. \eqref{ctomc}.

\section{Models for graphs with motifs}\label{rgmotifs}
In this section we focus on the case where the we have atoms of order larger than 2 which can be used as models for networks with extensive numbers of triangles and other bi-connected subgraphs. For this we place subgraph configuration ensembles within the context of some existing random graph models and representations \cite{Karrer2010RandomSubgraphs,Bollobas2011SparseClustering}. Establishing this connection allows various methods developed in the context of such models to be carried over to subgraph configuration models. 

\subsection{Kernel models}
In Ref. \cite{Bollobas2011SparseClustering} Bollob\'as, Janson and Riordan introduce a class of random graph models that generalizes non-homogeneous random graphs \cite{Bollobas2007TheGraphs} to the case where not only edges but also copies of small atomic subgraphs are added on to the vertices of the graph during the generation process. 
Given a set of atoms $M$ and a feature space $S$ a graph is generated by adding every possible embedding of $m \in M$ where vertex $i$ is mapped onto $v_i$ to the graph independently with probability $p_m$:
\begin{equation}
    p_m(v_1,v_2....,v_{|m|})= \frac
    {\mathcal{K}_m(s(v_1),s(v_2),...,s(v_{|m|}))}{N^{1-|m|}},
\end{equation}
where $s(v_i)\in S$ are vertex features and $\mathcal{K}$ is a function from $S^{|m|}$ to $[0,\infty)$. The normalization by $N^{1-|m|}$ ensures the graphs are sparse i.e. have on average $O(N)$ edges. 

Note that the subgraph configuration formulation slightly deviates from the one given above \cite{Bollobas2011SparseClustering} which considers embeddings of atomic subgraphs that is one to one mappings of vertices of atoms to the vertices of the graphs. This results in every $m$-subgraph to be considered $|\mathrm{Aut(m)}|$ times for addition. However, the formulations are essentially equivalent and a subgraph configuration model can be obtained by simply combining all embeddings that correspond to the same subgraph. Conversely, any kernel model over subgraph configurations with a bounded kernel $\mathcal{K}_m$ can be mapped onto a equivalent embedding based model by simply dividing $\mathcal{K}_m$ by $|\mathrm{Aut(m)}|$. A more detailed discussion can be found in \cite{Bollobas2011SparseClustering}. 

Canonical ensembles fall into the category of kernel models.  For the degree corrected model the vertex features are the expected atomic degrees $k_{m,i}(v)$ and the kernel is given by $\mathcal{K}_m=\frac{\sum_i \langle k_{m,i} \rangle}{|m|\prod_{i}\langle k_{m,i}\rangle}\prod_{v}k_{m,O(v)}(v)$ whereas for the homogeneous model the kernel is the constant and equal to $\langle k_m \rangle/|m|$ where $\langle k_m \rangle=n_m/N$. In \cite{Bollobas2011SparseClustering} the authors derive extensive results for properties of kernel models including component sizes, percolation properties, the degree distribution and subgraph counts. These results carry over to canonical subgraph configuration ensembles with little modification. 

\subsection{Microcanonical models}

Here we present some variations of the microcanonical models (Sec. \ref{Mcanonical}) which correspond to relaxations of constraints on the atomic degree distribution. These models in general require less parameters for the same set atoms compared to the model that conserves the atomic degree distribution at the level of orbits. First we consider the microcanonical ensemble where all the orbits that correspond to the same atom are aggregated:

\begin{equation}
    d_m(v)=\sum_{i}d_{m,i}(v),
\end{equation}
for all $m\in M$. Note that this is equivalent to removing the distinction between orbits of the same atom. Applying the corresponding transformations derived in Sec. \ref{mcagg} we obtain the following expression for the entropy:
\begin{equation}\label{smd}
\begin{split}
& S(M,\mathbf{d}_m)
=\sum_m\Bigg[-\ln(n_m!)-n_m\ln(|\mathrm{Aut(m)}|)\\&+\ln((|m|n_m)!)-\sum_v \ln(d_{m}(v)!)\\&
-\frac{|\mathrm{Aut(m)}|n_m^2}{2(n_m|m|)^{|m|}}\Big(\frac{\langle d_{m}^2\rangle}{\langle d_{m}\rangle}-1\Big)^{|m|}
-\frac{|m|-1}{2}\Big(\frac{\langle d_{m}^2\rangle}{\langle d_{m}\rangle}-1\Big) \Bigg],
\end{split}
\end{equation} 
where we used the fact that $\sum_i|O_{m,i}|=|m|$. Note that this is the same expression one would obtain if all atoms had a single orbit.

The model can be further relaxed so that only the total number of atoms attached to each vertex is conserved by aggregating atomic degrees of all orbits:
\begin{equation}
    d(v)=\sum_{m\in M}\sum_i d_{m,i}(v).
\end{equation}
Again using the transformations derived in Sec. \ref{mcagg} we obtain the following expression for the entropy:
\begin{equation}\label{std}
\begin{split}
 S(M,\mathbf{n}_m,&\mathbf{d})
=\ln((\sum_m|m|n_m)!)-\sum_v \ln(d(v)!) \\ &-\sum_m\Bigg[\ln(n_m!)+n_m\ln(|\mathrm{Aut(m)}|)
\\&+\frac{|\mathrm{Aut(m)}|n_m^2}{2(\sum_{m'}|m'|n_{m'})^{|m|}}\Big(\frac{\langle d^2\rangle}{\langle d\rangle}-1\Big)^{|m|}\\&
+\frac{|m|-1}{2}\frac{|m|n_m}{(\sum_{m'}|m'|n_{m'})}\Big(\frac{\langle d^2\rangle}{\langle d\rangle}-1\Big) \Bigg].\end{split}\end{equation} 
This model has only a single degree sequence as its parameter and hence the expression above and Eq. \eqref{Cent}, Eq. \eqref{smd}  all become equivalent when the model has a single atom with one orbit. 

Properties of microcanonical models including subgraph counts, component sizes and percolation properties can be found in Ref.\cite{Karrer2010RandomSubgraphs} which uses a generalization of the  generating function formalism for the edge configuration model. For more recent results spectral properties of the microcanonical model see Ref. \cite{Newman2019SpectraLoops}. These methods and results can be carried over to all the variants of the microcanonical model with minor modifications.

\subsection{Simplicial complexes, hypergraphs and bipartite models}
We now consider some widely used models that include higher order interactions in the form of cliques. First we consider simplicial complexes. Statistical ensembles of simplicial complexes have been studied before for example in Ref. \cite{Courtney2016GeneralizedComplexes}.

Simplicial complexes consisting of $d$ dimensional simplices are equivalent to subgraph configurations consisting of cliques of size $d+1$. Hence, such models can be recovered by considering  models for which $M$ consists of only $K_{d+1}$ \cite{Courtney2016GeneralizedComplexes}. Similarly, various hypergraph ensembles can be obtained by considering atoms that are cliques. 

Bipartite models are another type of model with atoms that consist of cliques. Bipartite representations have traditionally been used as models for collaboration networks \cite{Newman2001RandomApplications}. In the bipartite representation one has two sets of vertices one representing the authors and the other representing scientific publications. An edge between $i$ and $j$ indicates that $i$ is an author of $j$. The collaboration network between authors is obtained by projecting the bipartite representation on to the set of authors by connecting all authors that have co-authored a publication. The bipartite model consists of randomizing a given bipartite representation such that the degrees of the vertices in both partitions are conserved.  

The entropy of the bipartite model can be obtained by a subgraph configuration model where we have a single atom consisting of a vertex labelled single edge with its vertices having two distinct labels $t$(top) and $b$ (bottom). This atom has $|\mathrm{Aut(m)}|=1$ and two orbits of size 1. Using this Eq. \eqref{Cent} we obtain:
\begin{equation}\label{Sbi}
\begin{split}
S(\mathbf{d}_t,\mathbf{d}_b)&= \ln(n_e!) - \sum_{v\in T} \ln(d_t(v)!) -\sum_{v\in B}\ln(d_b(v)!) \\&-\frac{1}{2}\Big(\frac{\langle d_t^2 \rangle}{\langle d_t \rangle}-1\Big)\Big(\frac{\langle d_b^2 \rangle}{\langle d_b \rangle}-1\Big),
\end{split}
\end{equation}
where $n_e$ is the number of edges and $d_t$ and $d_b$ the degrees of vertices in the top and bottom partitions. Note that the self loop term vanishes as expected in the bipartite case. 

Bipartite representations are equivalent to subgraph configurations consisting of cliques where every $n$-clique corresponds to a $n$-author publication. Single author papers can be included in the model in the form of self loops with a single vertex. Conserving the degree of the top vertices is equivalent to conserving the number of publications for authors and conserving the degrees of bottom vertices is equivalent to conserving the number of $n$-author publications in the model. Consequently, bipartite models are equivalent to clique configuration models where all orbit degrees are combined. However, one important feature of the bipartite model is that the publications are assumed to be distinguishable. As the set of atoms $M$ is the set of publications in general it will contain multiple distinguishable cliques of the same size. Moreover every atom occurs once in the model. Cliques have one orbit and $|\mathrm{Aut(m)}|=|m|!$. Applying Eq. \eqref{std} we obtain:
\begin{equation}\label{Sbiclq}
\begin{split}
    S(M,\mathbf{d})&=\ln((\sum_{m}|m|)!)-\sum_v\ln(d(v)!)-\sum_m\ln(|m|!)\\&-\frac{1}{2}\sum_m\frac{|m|^2-|m|}{\sum_{m'}|m'|}\Big(\frac{\langle d^2\rangle}{\langle d\rangle}-1\Big).
\end{split}
\end{equation}
Note that in this case the correction term for multiple subgraphs vanishes since the model only contains a single copy of every atom. Although Eq.\ref{Sbi} concerns edge configurations in a bipartite graph and Eq. \eqref{Sbiclq} clique configurations the two equations are identical given that $\mathbf{d}_t=\mathbf{d}$ and $M$ consists of cliques of which the sizes are given by $\mathbf{d}_b$. 

Bipartite models also been advocated as general models of complex networks that have high clustering \cite{Guillaume2004BipartiteNetworks}. In most such cases the bipartite representation is not known in advance and a bipartite representation has to be inferred from the network instead. In this case it might be more suitable to assume that cliques of the same size are indistinguishable which would result in an additional term ($-\sum_m \ln(n_m!)$) in the entropy.

\subsection{Directed hypergraphs and power graphs}
There exist several alternative definitions for directed hypergraphs \cite{Gallo1993DirectedApplications,Ausiello2017DirectedSurvey,Jost2019HypergraphNetworks}. In general though directed hyperedges can be represented as directed atoms that consist of two sets of vertices such that all the vertices in one partition are connected to all vertices in the other partition via directed edges and alternative definitions differ in regards to the type of directed hyperedges they allow. Consequently, when set up with such atoms, subgraph configuration ensembles are equivalent to ensembles of directed hypergraphs. 

Another class of graph representations related to subgraph configurations are power graphs \cite{Ahnert2015GeneralisedNetworks} which represent networks as collections of cliques and complete bipartite subgraphs. Although initially not conceived as generative models it is possible to model graphs that have various types of power graph representations using subgraph configuration models that contain only cliques and bipartite cliques. 

\subsection{Statistical inference}
Although in principle subgraph configurations can model a large variety of higher order graph structures in many applications information on higher order interactions is not readily available and has to be inferred from pairwise interactions i.e. a graph. Even when when data on higher order interactions is available as in the case of bipartite representations scientific collaboration networks the data might not contain all forms of higher order interactions and it might be possible to infer these for the data. For instance, scientific collaboration networks might well contain higher order interaction patterns beyond cliques.  

Even in the setting where the atoms are known in advance inferring higher order interactions from graph data is a non-trivial problem. Although there exist some heuristics for extracting higher order interactions from graphs in the case of bipartite representations \cite{Guillaume2004BipartiteNetworks} and simplicical complexes \cite{Kahle2009TopologyComplexes} these are restricted to clique like interactions in undirected networks. An alternative and arguably more principled approach to obtaining higher order interactions in networks is to use statistical inference. In general inferring a subgraph configuration for a given network also involves finding the set atoms that is most appropriate for representing the given network. In this context the explosion of potential atoms as the order of subgraphs increases poses theoretical and computational challenges. Nevertheless when combined with non-parametric priors similar to those used in inference based methods in community detection \cite{Peixoto2017NonparametricModel} it is possible to perform Bayesian inference for subgraph configurations. Such inference procedures for atomic substructures and subgraph configurations based on the presented models is beyond the scope of this article and will presented in a separate publication.

\section{Networks with community structure}\label{SBM}
In this section we discuss the case where atoms have vertex and edge labels. At first we consider only single edge atoms and obtain models that include the SBM and many of its variations. Later we also discuss the case of higher order atoms and its potential implications for community detection methods based on statistical inference of SBMs. 

\subsection{SBMs and edge atoms with vertex labels}
In the presence of vertex and edge labels the definition of graph isomorphisms should be modified to include the preservation of vertex and edge labels. We first consider the case where the model only contains 2-edges. We denote single the edge atom with vertex labels $r$ and $s$ as $e_{rs}$. For $r=s$ we have $|Aut(e_{rr})|=2$ and a single orbit $O_{rr,r}$ of size 2 while for $r\neq s$ we have  $|Aut(e_{rs})|=1$ and two orbits each corresponding to a vertex label which we denote as $O_{rs,r}$ and $O_{rs,s}$.

\subsubsection{The homogeneous SBM} 
The SBM \cite{Holland1971TransitivityGroups.} is the standard model for networks with community structure. In the SBM vertices are assigned to one of $B$ blocks and edges are i.i.d  with success probability $p_{rs}$ for edges connecting two vertices from blocks $r$ and $s$. It can be shown that the SBM is the maximum entropy ensemble given block assignments of vertices and the expected number of edges between blocks \cite{Peixoto2012EntropyEnsembles}.

\subsubsection{Degree-corrected SBM}
The degree-corrected SBM can be obtained by assuming that every vertex has a unique label indicating its block membership and that the vertex labels of edges have to match block labels of the vertices.  
In the orbit degree corrected model fixing the atomic degree of a vertex $v$ in block $r$ $\mathbf{d}(v)=(d_1,...,d_B)$ is equivalent to giving the number of neighbours of different types $v$ has. Under these specifications the entropy of the degree-corrected SBM can be obtained by Eq. \eqref{Cent}.

However, in most applications rather then specifying community specific degrees, which would require $B$ degrees for each vertex, the total degree of vertices is specified instead. The entropy for such models can be obtained by aggregating the degrees of all orbits corresponding the same vertex label resulting in the following constraints:
\begin{equation}
    d_r(v)=\sum_{s} d_{rs,r}(v),
\end{equation}
for vertices $v$ in block $r$. Under these constraints we get the following expression for the microcanonical entropy:
\begin{equation}\label{SDCSBM}
\begin{split}
    &S(\{\mathbf{d}_r\},\{n_{rs}\})=-\sum_{rs}[\ln(n_{rs}!)+n_{rs}\ln(|Aut(e_{rs})|)]\\ &+\sum_r\ln((\sum_s|O_{rs,s}|n_{rs})!)-\sum_{r,v}log(d_r(v)!)\\
    &-\sum_{rs}\frac{\mathrm{|Aut(e_{rs})|}n_{rs}^2}{2}\prod_{i}\frac{1}{n_i^{|O_{rs,i}|}}\Big(\frac{\langle d_{i}^2\rangle}{\langle d_{i}\rangle}-1\Big)^{|O_{rs,i}|}\\
    &-\sum_r\frac{n_{rr}}{\sum_s|O_{rs,s}|n_{rs}}(\frac{\langle d_{r}^2\rangle}{\langle d_{r}\rangle}-1),
\end{split}
\end{equation}
where $n_{rs}$ is the number of edges between blocks $r$ and $s$ and $n_i=\sum_t n_{it}|O_{it,i}|$ is the number of half edges with vertex label $i$. Substituting $|O_{rs,s}|=1+\delta_{rs}$ and $\mathrm{|Aut(e_{rs})|}=1+\delta_{rs}$ with their numerical values and applying Stirling's approximation to the factorial terms we recover the expression derived in Ref. \cite{Peixoto2012EntropyEnsembles}. The entropy for the canonical ensemble can also be obtained in a similar fashion and agrees with known expressions \cite{Peixoto2012EntropyEnsembles}.

It is also possible to formulate SBMs with intermediate parametric complexity. For instance, one could construct a SBM where one distinguishes between in community degree and out community degree by aggregating only the orbits of the edge atoms that have two distinct vertex labels. Consequently the entropy of such a model is the sum of two copies of Eq.\ref{SDCSBM} one where $n_{rs}$ are set to 0 for $r=s$ plus one where $n_{rs}$ are set to 0 for $r\neq s$. 

\subsubsection{Overlapping SBM}
The SBM with overlapping blocks \cite{Peixoto2015InferringNetworks} can be obtained by relaxing the condition that each vertex can only receive orbits that correspond to single block label. In other words the atomic degree vector is allowed to have non-zero entries for multiple block labels. As such the non-overlapping SBM is a subset of the overlapping SBM and the distinction between the two models arises only due to additional assumptions about the atomic degree sequence and hence the entropy expressions are identical for both variants of the SBM \cite{Peixoto2015InferringNetworks}. 

\subsubsection{Directed SBM}
The directed SBM differs from the undirected SBM only with respect the number of atoms involved and their symmetries. In the directed case there are two types of directed edges for every pair of distinct labels and all directed single edge atoms have 2 orbits and their automorphism groups are trivial. In the directed case it is customary to conserve the in and out degrees of vertices separately which is equivalent to placing orbits corresponding to the same vertex label with an incoming and outgoing edge into separate groups. Making these changes the entropy of the directed SBM \cite{Peixoto2012EntropyEnsembles} can be recovered following the same procedure as in the undirected case.

\subsection{Edge labels: link communities/hidden layers}
We consider two cases of edge labels, the first being the case where the model produces an unlabelled graph and edge labels are hidden variables to be inferred from the data similar to vertex labels in the SBM. This essentially provides the counterpart of the SBM for link communities \cite{Evans2009LineCommunities,Ahn2010LinkNetworks}. The second case which we shall consider later is multi-layered networks where edge labels correspond to different layers in a multi-layer networks. Therefore the model with labelled edges can also be interpreted as one with hidden layers \cite{Valles-Catala2016MultilayerNetworks}. 

In the case of labelled edges one simply obtains a model with independent layers where each layer is a (edge) configuration model. As a result the entropy of such a model is simply the sum of the entropies of these  models. Note that in this case the community/layer membership of a vertex can be deduced from its atomic degrees. 

The link community model is equivalent to the SBM with overlapping communities where the counts of the atoms corresponding to inter community edges is set to zero. This is due to the fact that the edge labelled single edge has the same automorphism group as the vertex labelled single edge atom that has the same label on both of its vertices. Note that in link community models communities can still be connected either through direct overlap or via other link communities. Although link community models are a subset of the overlapping SBMs link community models in general require less parameters than the overlapping SBM with the same number of blocks.  

\subsubsection{Edge atoms with vertex and edge labels}\label{NElabel}

In the  degree corrected case the link community the model is equivalent to having multiple independent (edge) configuration models that are coupled through vertex intersections. Hence if every edge label is assumed to have its own exclusive set of vertex labels the model becomes equivalent to the case where each link community consist of a SBM which includes the case of bipartite link communities. On the other hand if vertex labels are shared across edge labels one obtains a SBM with multiple layers that have a common block structure. However these do not exhaust the possible models as in principle the model does not intrinsically restrict the relation between vertex and edge labels.  

\subsection{Model selection}

Considering models with vertex and edge labelled atoms leads to variety of generative models for networks with communities corresponding to different notions of network communities. Formulating them within an unified framework enables principled model selection as different variants of the model share the same type of parameters. As has been for instance done in the case of overlapping and non-overlapping SBMs in Ref.\cite{Peixoto2015ModelGroups} using Bayesian non-parametrics and by setting up all model variants with the same priors. 

\subsection{Higher order atoms and network communities}
Including higher order structures such as triangles in generative models can have a significant impact on inference of network communities. In general atomic structures whether known a priory or inferred from data provide additional information that can either reinforce or counteract evidence for community structures. For instance consider a random graph with $N$ vertices onto which we add $N/6$ triangles and $kN$ edges at random resulting at most $N/2$ vertices that have a triangles attached to them. Consequently, the difference in density between the vertices having triangles attached to them and those which have none might be misinterpreted as evidence for the presence community structure. 
On the other hand higher order structures might also reinforce and facilitate the detection of communities when the distribution of atoms is strongly correlated with block structures and hence facilitate the detection of block structures in data that otherwise might not be detectable by SBMs that only consider pairwise interactions. For instance the method proposed in Ref. \cite{Benson2016Higher-orderNetworks.} leverages subgraph structures for the purpose community detection. Though these methods assume that the topology of the hyperedges is known in advance and that every such subgraph constitutes a hyperedge. The presented models can be used to infer an optimal (generalized) hypergraph representation for a given network without requiring any prior knowledge on the topology of atoms/hyperedges which then could be further used as a basis for community detection.

\subsubsection{Hypergraph communities}
Detecting communities in hypergraphs is an active area of research \cite{Papa2007HypergraphClustering,Kaminski2019ClusteringModularity.} and the presented approach can be used to generalize the SBM by considering hyperedges with vertex and/or edge labels. This could be used to generalize inference based  community detection methods  \cite{Peixoto2017NonparametricModel} to hypergraphs and directed hypergraphs. 

\section{Multilayer networks}\label{Mlayer}
Multilayer networks \cite{Bianconi2018MultilayerFunction} can be modelled as networks with labelled edges where edge labels indicate layer membership. In subgraph configuration models for multilayer networks the projection from configurations to multilayer graphs hence should preserve edge labels that correspond to layer assignments. We shall mostly focus on cases where layers are correlated since the uncorrelated case layers are independent and the entropy can be obtained by summing over the entropy of each layer specific model.  

\subsection{Vertex couplings of layers}
One way of obtaining correlated layers is to couple the degree distribution across layers. Groups of layers can be coupled by fixing the degrees of vertices in the graph obtained by aggregating the layers in group under consideration. Given such a group of layers $L$ this is equivalent to considering constraints of the form $\{e_l=n_l\}$ and $d_L(v)=\sum_{l\in L} d_l(v)$ where $e_l$  is the total number of links and $d_l$ is the degree in layer $l$. 

In the canonical ensemble layers are independent given the expected aggregate degrees of vertices and the entropy has the from of layer specific models with expected degrees given by $d_l(v)=\frac{n_l}{n_L}d_L(v)$ where $n_L=\sum_{l} {n_l}$. In the microcanonical ensemble the entropy is the sum of unlabelled (edge) configuration model with the aggregated degree distribution $\mathbf{d}_L$ plus $\ln({n_L\choose\mathbf{n}_l })$ which accounts for the assignment of edges to layers in the right proportions. Hence network generation and layer assignment of edges are independent in the microcanonical model. Such models have previously been considered in the context of multi-layer SBMs \cite{Peixoto2015InferringNetworks}. As argued in Ref. \cite{Peixoto2015InferringNetworks} in such models layers can be interpreted as edge features of the aggregated graph.  

\subsection{Edge couplings of layers}

Layers can be coupled at the level of edges by considering atoms that consist of multiple parallel edges from different layers. This is essentially equivalent to creating a new layer consisting of a specific multi edge pattern of layers and the entropy of such models is given by the sum over individual entropies corresponding to the label combinations/atoms in the model. 

Maximum entropy models for constraints on multi layer intersection patterns have been previously studied in Ref. \cite{Bianconi2013StatisticalOverlap}. In Ref. \cite{Bianconi2013StatisticalOverlap} the ensemble is constrained in terms of counts of (induced) copies of multi layer intersection patterns resulting in a model where pairs of vertices are connected by a unique multi edge pattern. As a result this model necessarily includes all intersection patterns found in a given multilayer graph. Subgraph configuration models on the other hand allow for intersection patterns not explicitly included in the model to be formed by random intersections of other atoms. As a result subgraph configuration models could be used to infer a more concise subset inter-layer edge intersection patterns for a given graph. Nevertheless, if both models are set up with the same set of intersection patterns one recovers the entropy expression derived in Ref. \cite{Bianconi2013StatisticalOverlap}. Models that combine both vertex and edge couplings can be obtained by aggregating the atomic degrees in models with edge couplings.

\subsection{Higher order couplings}
It is possible to construct subgraph configuration models where layers are coupled beyond vertex and edge couplings. For instance, two layers can be coupled by considering a triangular atom where two edges have one label and the third edge another label. In multilayer networks the number of potential atoms can be very large even for small atoms especially for networks with many layers \cite{Battiston2017MultilayerNetworks}. Hence, having models that can characterize the local structures in such networks concisely are of great interest. 

\subsection{Multilayer SBMs}
As in the single layer case different variants of the multilayer SBM can be obtained by considering single atoms with both edge and vertex labels. Models where group labels are shared across all layers were considered in Ref.s \cite{Peixoto2015InferringNetworks,Bianconi2013StatisticalOverlap}. Such models fall into the category of models with edge atoms that have both vertex and edge labels which were discussed in Sec.\ref{NElabel}.

\section{Conclusions}\label{conc}

We have formulated a general class of maximum entropy models for higher order network interactions. This class is based on explicitly representing network interactions by atomic subgraphs which can have a large variety of topological features. The resulting models are analytically tractable and can generate networks with nontrivial and complex subgraph structures. We calculated general expressions for the entropy of such models and presented a coarse graining procedure based on aggregating ensemble constraints that allows the parametric complexity of the models to be controlled.

We have shown that the presented models include a large variety of models from the literature ranging from models for networks with motifs to SBMs to various models of multilayer networks and that the entropy of all these models can be recovered from a single expression that is characterized by the symmetry groups of atomic substructures. We also identified new generative models for structural features such as link communities and multilayer network motifs, this showing the generality of our methodology.

The presented models thus offer a powerful and flexible framework where network structures are represented in terms of atomic substructures. Subgraph configuration models also allow for the construction of models that combine multiple network features that better reflect empirical features of complex real-world networks.

One of the main applications of the presented models is in the area of inference based methods similar to those developed in the context of SBMs \cite{Peixoto2013ParsimoniousNetworks}. The presented models significantly expand the types of network structures that can be addressed by statistical inference methods. Moreover having a general class of models that includes multiple alternative models provides a consistent framework for comparing alternative models and representations via model selection. Such methods have been used to the discriminate between overlapping and non-overlapping communities \cite{Peixoto2015ModelGroups} and alternative models for communities in multilayer networks \cite{Peixoto2015InferringNetworks}. In the context of subgraph configuration models the objective of inference is to identify a set of atoms together with a set of constraints on their distribution that optimally describe the observed data. Hence, the development such inference algorithms remains an important area for future research.

\section*{Acknowledgements}
Anatol E. Wegner and Sofia Olhede acknowledge the financial support of  the UK's Food Standards Agency. Sofia Olhede acknowledges support fro the European Research Council under Grant CoG 2015-682172NETS within the Seventh European Union Framework Program. 


\bibliography{main}

\end{document}